\documentclass[12pt]{amsart}

\newtheorem{fact}{Fact}[section]
\newtheorem{lemma}[fact]{Lemma}
\newtheorem{theorem}[fact]{Theorem}
\newtheorem{definition}[fact]{Definition}
\newtheorem{example}[fact]{Example}
\newtheorem{remark}[fact]{Remark}

\newcommand{\CC}{{\mathbb C}}
\DeclareMathOperator{\N}{\mathbb N}
\DeclareMathOperator{\C}{\mathcal C }
\DeclareMathOperator{\A}{\mathcal A}
\DeclareMathOperator{\n}{\mathfrak{n}_{--}}
\DeclareMathOperator{\np}{\mathfrak{n}_{+}}
\DeclareMathOperator{\g}{\mathfrak{g}}
\DeclareMathOperator{\F}{\mathcal F}
\DeclareMathOperator{\h}{\mathfrak{h}}
\DeclareMathOperator{\res}{Res} \DeclareMathOperator{\sym}{Sym}
\DeclareMathOperator{\id}{id} \DeclareMathOperator{\asym}{ASym}
\DeclareMathOperator{\sgn}{sgn} \DeclareMathOperator{\dlog}{dlog}

\begin{document}

\begin{abstract} We study the canonical $U(\n)$-valued differential form, whose projections
to different Kac-Moody algebras are key ingredients of the
hypergeometric integral solutions of KZ-type differential equations
and Bethe ansatz constructions. We explicitly determine the
coefficients of the projections in the simple Lie algebras $A_r,
B_r, C_r, D_r$ in a conveniently chosen Poincar\'e-Birkhoff-Witt
basis.
\end{abstract}

\title[ PBW Expansions of the Canonical
Differential Form]{Combinatorics of Rational Functions and
Poincar\'e-Birkhoff-Witt Expansions of the Canonical $U(\n)$-valued
Differential Form}
\author{R. Rim\'anyi}
\address{Department of Mathematics, University of North Carolina at Chapel Hill}
\email{rimanyi@email.unc.edu}
\author{L. Stevens}
\address{Department of Mathematics, Universita' di Roma "La
Sapienza"}
\email{stevens@mat.uniroma1.it}
\author{A. Varchenko}
\address{Department of Mathematics, University of North Carolina at Chapel Hill}
\email{anv@email.unc.edu}

\thanks{\noindent Supported by
        NSF grant DMS-0405723 (1st author), DMS-0244579 (3rd author) \\
Keywords:
 canonical differential form,
KZ equation, Bethe ansatz, PBW-expansion, symmetric rational functions\\
AMS Subject classification 33C67}

\maketitle

\section{Introduction}\label{intr}

For a Kac-Moody algebra $\g$, let $V$ be the tensor product
$V_{\Lambda_1}\otimes \ldots\otimes V_{\Lambda_n}$ of highest weight
$\g$-modules. The $V$-valued hypergeometric solutions of
Knizhnik-Zamolodchikov-type differential equations have the form
\cite{sv1}, \cite{sv}:

\begin{equation}\label{integral}
I(z)\ =\ \int_{\gamma(z)} \ \Phi(t,z)\ \Omega^{V}(t,z)\ .
\end{equation}
Here $t=(t_1,\ldots,t_k)$, $z=(z_1,\ldots,z_n)$, $\Phi$ is a
scalar multi-valued (master)  function, $\gamma(z)$ is a suitable cycle in
$t$-space depending on $z$, and $\Omega^{V}$ is a $V$-valued rational
differential $k$-form.

The same $\Phi$ and $\Omega^V$ have applications to the Bethe Ansatz method.
It is known \cite{asymp} that the values of $\Omega^V$ at the
critical points of $\Phi$ (with respect to $t$) give eigenvectors of
the Hamiltonians of the Gaudin model associated with $V$.

For every
$V =V_{\Lambda_1}\otimes \ldots\otimes V_{\Lambda_n}$,
the $V$-valued differential form
$\Omega^V$ is
constructed out of a single $U(\n)$-valued
differential form $\Omega^{\g}$, where $\g=\n \oplus \h \oplus \np$ is
the Cartan decomposition of $\g$, and $U(\n)$ denotes the universal
enveloping algebra of the Lie algebra $\n$, see Appendix. In applications, it is important to have convenient
formulas for $\Omega^{\g}$, and this is the goal of the present paper.

In \cite{matsuo}, Matsuo
 suggested a formula $\int
\Phi(t,z)\widetilde{\Omega}^{V}(t,z)$ for solutions of the KZ
equations for  $\g=sl_{r+1}$. His differential form
$\widetilde{\Omega}^V$ also can be
constructed from a $U(\n)$-valued form $\widetilde{\Omega}^{sl_{r+1}}$
in the same way  as $\Omega^V$ from $\Omega^{\g}$.
It is known that for $sl_2$, the two forms
\begin{align}\notag
\Omega^{sl_2} & =\sum_{k=0}^{\infty}\Big(\sum_{\pi\in\Sigma_k} \sgn(\pi)\cdot
\frac{dt_{\pi(1)}}{t_{\pi(1)}}\wedge\frac{d(t_{\pi(2)}-t_{\pi(1)})}{t_{\pi(2)}-t_{\pi(1)}}\wedge\ldots
\wedge\frac{d(t_{\pi(k)}-t_{\pi(k-1)})}{t_{\pi(k)}-t_{\pi(k-1)}}\Big) \otimes f^k,
\\ \notag
\widetilde{\Omega}^{sl_2} & =\sum_{k=0}^{\infty} \Big(\bigwedge_{i=1}^k \frac{dt_i}{t_i}\Big)\otimes f^k
\end{align}
coincide. For  $r>1$, the form $\Omega^{sl_{r+1}}$ is a polynomial
in $f_{1},\dots,f_{r}$ with scalar differential forms as
coefficients, while the Matsuo form $\widetilde{\Omega}^{sl_{r+1}}$
is a sum over a Poincar\'e-Birkhoff-Witt basis of $U(\n)$ with
coefficients of the same type. Both forms have some advantages.  The
form $\Omega^{\g}$ is given by the same formula for any $\g$. The
formula for $\widetilde{\Omega}^{sl_{r+1}}$ has less terms and less
apparent poles (see above for the apparent poles at $t_i-t_j=0$).
The advantage of having an expression in terms of a PBW basis is
most spectacular for representations with 1-dimensional
weight-subspaces.

In this paper, we prove that
$\widetilde{\Omega}^{sl_{r+1}}=\Omega^{sl_{r+1}}$ and give similar
Poincar\'e-Birkhoff-Witt expansions for the differential form
$\Omega^{\g}$  for the simple Lie algebras $\g$ of types $B_r$,
$C_r$, $D_r$.

As a byproduct, we obtain results on the combinatorics of rational functions. Namely, some
non-trivial identities are established among certain rational functions with partial symmetries.
The results are far reaching generalizations of the prototype of these formulas, the
``Jacobi-identity''
$$\frac{1}{(x-y)(x-z)}+\frac{1}{(y-x)(y-z)}+\frac{1}{(z-x)(z-y)}=0.$$

In all $A_r, B_r, C_r, D_r$ cases, the coefficients of
$\Omega^{\g}$ can be encoded by diagrams relevant to sub-diagrams
of the Dynkin diagram of $\g$. One may expect that the same phenomenon occurs in
a more general Kac-Moody setting, too.

\medskip

According to the formulas for $\Phi$ and $\Omega^V$ in \cite{sv1},
\cite{sv},
the
poles of $\Omega^V$ contain the singularities of $\Phi$. From our PBW
expansion formulas, it follows that the poles of $\Omega^V$ {\em
coincide} with the singularities of $\Phi$, hence it makes sense to
consider the values of $\Omega^V$ at (e.g.) the critical points of
$\Phi$, as is needed in
the
Bethe ansatz applications.

It was shown in \cite{trig} that the Matsuo type hypergeometric
solutions of the $sl_{r+1}$ KZ-equations satisfy the complementary
dynamical difference equations. According to our result
$\Omega^{sl_{r+1}}=\widetilde{\Omega}^{sl_{r+1}}$
,
the hypergeometric
solutions (\ref{integral}) also satisfy the dynamical difference
equation for $\g=sl_{r+1}$.

In \cite{fv}, the hypergeometric solutions $I(z,\lambda)=\int
\Phi^{ell}(t,z)\Omega^{V,ell}(t,z,\lambda)$ of the KZB equations
were constructed. Here $\Phi^{ell}$ is the elliptic scalar master
function depending on the same variables $t,z$ as $\Phi$, and
$\Omega^{V,ell}(t,z,\lambda)$ is the elliptic analogue of
$\Omega^V$, which is a $V$-valued differential form depending also
on $\lambda\in\h$. It would be useful to find PBW type expansions of
$\Omega^{V,ell}$ similar to our PBW expansions of $\Omega^V$.

The hypergeometric solutions of qKZ, the quantum version of KZ equations,
for $sl_{r+1}$-modules $V$ were described in \cite{qgeometry},
\cite{qsl} as $I(z)=\int\Phi_q(t,z) \Omega_q^V(t,z)$. There the
$V$-valued differential form $\Omega_q^V$ was given in a PBW
expansion. Our PBW formulas for the B, C, D series may suggest
integral formulas for solutions of B, C, D type qKZ equations.

The structure of the cycles $\gamma$ in (1) for arbitrary $\g$ was
analyzed in \cite{varbook}.  The cycles were presented as linear
combinations of multiple loops, and that presentation established
a connection between multi-loops and monomials $f_{i_k}\ldots f_{i_1}$
in $U_q(\n)$, where $U_q(\n)$ is the $\n$-part of the quantum group
$U_q(\g)$. That connection in particular
gives an identification of the monodromy of
the KZ equations with the $R$-matrix representations associated with
$U_q(\g)$. Our PBW expansions of $\Omega^{\g}$ suggest that there
might be an interesting PBW type geometric theory of cycles for each
$\g$, in which the cycles are presented by linear combinations of
cells corresponding to elements of the PBW basis in the corresponding
$U_q(\n)$.

It would also be interesting to compare our PBW formulas with
Cherednik's formulas for solutions of the trigonometric KZ equations
\cite{cher}.

\smallskip

The authors thank D. Cohen for helpful discussions.

\section{Symmetrizers, signs and other conventions}

\subsection{Symmetrizers} For a nonnegative integer $r$ and
$k=(k_1,\ldots,k_r)\in \N^r=\{0,1,2,\ldots\}^r$, we will often consider
various `objects' $x(t^{(i)}_j)$
(functions, differential forms, flags),
depending on the $r$ sets of variables
\begin{equation}\label{tij} t^{(1)}_1, t^{(1)}_2,\ldots,t^{(1)}_{k_1},\qquad
t^{(2)}_1, t^{(2)}_2,\ldots,t^{(2)}_{k_2},\qquad \ldots \qquad
t^{(r)}_1, t^{(r)}_2,\ldots,t^{(r)}_{k_r}.\end{equation}
Let $G_k$ be the product $\prod_i \Sigma_{k_i}$ of symmetric groups.
We define the action of $\pi\in G_k$ on $x$
by permuting the $t^{(i)}_j$'s with the same upper indices. Then we
define the symmetrizer and antisymmetrizer operators
$$\sym_k x(t^{(i)}_j)=\sum_{\pi\in G_k} \pi \cdot x,\qquad\qquad
\asym_kx(t^{(i)}_j)=\sum_{\pi\in G_k} \sgn(\pi)\ \pi \cdot x.$$

Let $|k|=\sum_ik_i$. For a function (`multi-index') $J:\{1,\ldots,|k|\}$ $\to$
$\{1,\ldots,r\}$ with $\#J^{-1}(i)=k_i$, let $c:\{1,\ldots,|k|\}\to
\N$ be the unique map whose restriction to $J^{-1}(i)$ is the
increasing function onto $\{1,\ldots, k_i\}$. Then $J$ defines an
identification of $(t_1,\ldots,t_{|k|})$ with the variables in (\ref{tij})
by identifying
\begin{equation}\label{id}  t_u  \qquad \hbox{with}\qquad t_{c(u)}^{(J(u))}.\end{equation}
Thus, if $x$ depends on
$t_1,\ldots, t_{|k|}$ and $J$ is given, we can consider $x$ depending on
the variables in (\ref{tij}).  For example, we can (anti)symmetrize $x$:
$$\sym^J_k x(t_u)=\sym_k x(t_j^{(i)}), \qquad\qquad \asym^J_k x(t_u)=\asym_k
x(t_j^{(i)}).$$

\subsection{The sign of a multi-index; volume forms} Let $J_0$ be the unique {\em increasing} function $\{1,\ldots,|k|\}\to
\{1,\ldots,r\}$ with $\#J_0^{-1}(i)=k_i$, and let $J$ be any function
$\{1,\ldots,|k|\}\to\{1,\ldots,r\}$ with $\#J^{-1}(i)=k_i$. Then the identifications defined in (\ref{id}) for $J$ and $J_0$
together define a permutation of $1,\ldots,|k|$. The sign of this permutation
will be denoted by $\sgn(J)$. E.g. $\sgn(1\mapsto 1, 2\mapsto 1, 3\mapsto 2)=1$,
$\sgn(1\mapsto 1, 2\mapsto 2, 3\mapsto 1)=-1$.

Define the `standard volume form' $dV_k$ to be $dt^{(1)}_1\wedge\ldots\wedge dt^{(1)}_{k_1}\wedge
dt^{(2)}_1\wedge\ldots\wedge dt^{(2)}_{k_2}\wedge\ldots\wedge dt^{(r)}_1\wedge\ldots\wedge dt^{(r)}_{k_r}$.
Observe that if we use the identification (\ref{id}), then $dt_1\wedge dt_2\wedge\ldots\wedge dt_{|k|}$ is equal
to sgn$(J)\cdot dV_k$.

\subsection{The star multiplication}\label{star}  For $k\in \N^r$, let $\F_k$ be the vector space of rational
functions in the variables in (\ref{tij}) which are symmetric under the action of $G_k$.
We define a multiplication (c.f. \cite[6.4.2]{cbms})
$*:\F_k\otimes \F_l \to \F_{k+l}$ by $$(f*g)(t^{(1)}_1,\ldots,t^{(1)}_{k_1+l_1},\ldots,t^{(r)}_1,\ldots,t^{(r)}_{k_r+l_r})=$$
$$\frac{1}{\prod_i k_i!l_i!}\sym_{k+l}\Big( f(t^{(1)}_1,\ldots,t^{(1)}_{k_1},\ldots,t^{(r)}_1,\ldots,t^{(r)}_{k_r}) \cdot
g(t^{(1)}_{k_1+1},\ldots,t^{(1)}_{k_1+l_1},\ldots,t^{(r)}_{k_r+1},\ldots,t^{(r)}_{k_r+l_r})\Big).$$
For example, if we write $t$ for $t^{(1)}$ and $s$ for $t^{(2)}$, then
$$\frac{1}{t_1t_2}*\frac{1}{t_1(s-t_1)}=
\frac{1}{t_1t_2}\cdot\frac{1}{t_3(s-t_3)}+\frac{1}{t_1t_3}\cdot\frac{1}{t_2(s-t_2)}+\frac{1}{t_2t_3}\cdot\frac{1}{t_1(s-t_1)}.$$
This multiplication makes $\oplus_k \F_k$ an associative and commutative algebra.

\section{Arrangements. The Orlik-Solomon algebra and its dual.
Discriminantal arrangements and their symmetries} \label{arr} Let
$\C$ be a hyperplane arrangement in $\CC^n$. In this section we
recall two algebraic descriptions of the cohomology of the
complement $U=\CC^n-\cup_{H\in\C}H$, as well as properties of the
discriminantal arrangement which will be needed later. The general reference is
\cite{sv}.

\subsection{The Orlik-Solomon algebra}
For $H\in \C$, let $\omega_H$ be the logarithmic differential form
$df_H/f_H$, where $f_H=0$ is a defining equation of $H$. Let
$\A=\A(\C)$ be the graded $\CC$-algebra with unit element
generated by all $\omega_H$'s, $H\in \C$. The elements of $\A$ are
closed forms on $U$, hence they determine cohomology classes.
According to Arnold and Brieskorn, the induced map $\A\to
H^*(U;\CC)$ is an isomorphism. The degree $p$ part of $\A$ will be
denoted by $\A^p$.

\subsection{Flags} Non-empty intersections of hyperplanes in $\C$ are called
edges. A $p$-flag of $\C$ is a chain of edges
$$F=[\CC^n=L^0\supset L^1\supset L^{2}\supset\ldots\supset L^{p-1}\supset
L^p],$$ where codim $L^i=i$. Consider the complex vector space
generated by all $p$-flags of $\C$ modulo the relations
$$\sum_{L}[L^0\supset\ldots\supset L^{i-1}\supset
L\supset L^{i+1}\supset\ldots\supset L^p]=0,\qquad (0<i<p),$$ where
the summation  runs over all
 codim $i$ edges $L$ that contain
$L^{i+1}$ and are contained in $L^{i-1}$. This vector space is
denoted by $Fl^p=Fl^p(\C)$, and let $Fl$ be the direct sum
$\oplus_p Fl^p$.

\subsection{Iterated residues} \label{itres}

According to \cite[Th.~2.4]{sv}, $\A$ and $Fl$ are dual graded vector spaces.
The value of a differential form on
a flag is given by an iterated residue operation
$$\res:Fl\otimes \A\to \CC,$$
defined as follows. Let
$F=[L^i]$ be a $p$-flag and $\omega\in \A$ a $p$-form on $U$.
Then
$$\res_{F}\omega=\res_{L^p}\Big(\res_{L^{p-1}}\big(\ldots
\res_{L^2}(\res_{L^1}(\omega))\ldots\big)\Big)\ \ \in\CC.$$

\subsection{The discriminantal arrangement and its symmetries}
The discriminantal arrangement $\C^n$ in $\CC^n$ is defined as the
collection of hyperplanes
$$t_i=0 \quad (i=1,\ldots, n)\qquad\qquad \text{and}\qquad\qquad
t_i-t_j=0 \quad (1\leq i<j \leq n).$$

Let us fix $r$ non-negative integers (weights)
$k=(k_1,\ldots,k_r)$ with $\sum k_i=|k|$ and consider $\CC^{|k|}$
with coordinates
$$(t^{(1)}_1,\ldots,t^{(1)}_{k_1},
t^{(2)}_1,\ldots,t^{(2)}_{k_2}, \ldots\ldots\ldots,
t^{(r)}_1,\ldots,t^{(r)}_{k_r}).$$ The group $G_k=\prod
\Sigma_{k_i}$ then acts on $\CC^{|k|}$ (by permuting the
coordinates with the same upper indices) which then induces an
action of $G_k$ on $\A(\C^{|k|})$ and $Fl(\C^{|k|})$.

The skew-invariant subspaces (i.e. the collection of $x$'s for
which $\pi\cdot x=\sgn(\pi)x \ \forall \pi\in G_k$) of
$\A^{|k|}(\C^{|k|})$ and $Fl^{|k|}(\C^{|k|})$ will be denoted by
$\A^{G_k}$ and $Fl^{G_k}$, respectively.
The duality stated in
\ref{itres} is consistent with the group-action in the sense that
$\A^{G_k}$ and $Fl^{G_k}$ are dual vector spaces.

\subsection{Flags of the discriminantal arrangement.} \label{flags}

Let $U_r$ be the free associative algebra generated by $r$ symbols
$\tilde{f}_1, \tilde{f}_2, \ldots, \tilde{f}_r$. It is multigraded
by $\N^r$; the $(k_1,\ldots,k_r)$-degree part will be denoted by
$U_r[k]=U_r[k_1,\ldots, k_r]$. For any non-zero homogeneous
element in $U_r[k]$, we define its {\em content} to be $k$. It is
proved in \cite[Th.~5.9]{sv} that $U_r[k]$ is isomorphic to
$Fl^{G_k}$ under the following map. For $J:\{1,\ldots,|k|\}\to
\{1,\ldots,r\}$ with $\#J^{-1}(i)=k_i$, the monomial
$\tilde{f}_J=\tilde{f}_{J(|k|)}\tilde{f}_{J(|k|-1)}\ldots
\tilde{f}_{J(2)} \tilde{f}_{J(1)}\in U_r[k_1,\ldots,k_r]$
corresponds to $\frac{\sgn(J)}{\prod_i k_i!}\asym_k^J(F)\in
Fl^{G_k}$, where $F$ is the $|k|$-flag
$$[\CC^{|k|}\supset (t_1=0) \supset (t_1=t_2=0)\supset
\ldots\supset (t_1=\ldots=t_{|k|-1}=0)\supset
(t_1=\ldots=t_{|k|}=0)]$$
with its variables $t_u$ identified with $t^{(i)}_j$'s
as defined by (\ref{id}).

\begin{example}\rm For $r=2$, $k=(1,1)$, we have the correspondence
$$\tilde{f}_2\tilde{f}_1 \leftrightarrow [\CC^2\supset (t_1^{(1)}=0) \supset
(t_1^{(1)}=t_1^{(2)}=0)],\qquad \tilde{f}_1\tilde{f}_2 \leftrightarrow
-[\CC^2\supset (t_1^{(2)}=0) \supset (t_1^{(2)}=t_1^{(1)}=0)].$$
For $r=2$, $k=(2,1)$, we have the correspondence
\begin{align} \notag
\tilde{f}_1^2\tilde{f}_2 \leftrightarrow   \frac{1}{2}  \Big( &
[\CC^3\supset(t_1^{(2)}=0)\supset (t_1^{(2)}=t^{(1)}_1=0)\supset
(t_1^{(2)}=t^{(1)}_1=t^{(1)}_2=0)] - \\\notag
& [\CC^3\supset(t_1^{(2)}=0)\supset (t_1^{(2)}=t^{(1)}_2=0)\supset
(t_1^{(2)}=t^{(1)}_2=t^{(1)}_1=0)]\Big).
\end{align}
\end{example}

\section{The canonical differential form}

Using the identifications of Section \ref{arr}, the tensor product
$$\A^{G_k}\otimes U_r[k_1,\ldots,k_r] $$
is the tensor product of a vector space with its dual space.
Therefore the canonical element, $\sum_i b_i^*\otimes b_i$ for any
basis $\{b_i\}$ of $U_{r}[k]$ and the dual basis $\{b_i^*\}$ of $\A^{G_k}$,
is well defined---it does not depend on the
choice of the basis of $U_{r}[k]$. We will call this element the {\em canonical
differential form of weight $k$} and denote it by $\Omega_k$.
Tracing back the identifications of Section \ref{arr}, we get the
explicit form.

\begin{theorem}\cite{sv} \label{canform} Let $t_0=0$. The canonical differential form is
\begin{align}\notag
\Omega_k& =  \sum_{J}  \sgn(J)\cdot \asym_k^J \Big(
  \bigwedge_{u=1}^{|k|} \dlog (t_u-t_{u-1}) \ \Big) \otimes \tilde{f}_J \\ \notag
 & =\sum_{J}  \sym_k^J \Big(\
  \prod_{u=1}^{|k|} \frac{1}{t_u-t_{u-1}} \ \Big) dV_k \otimes \tilde{f}_J\ \ \ \ \ \in \A^{G_k}\otimes U_r[k],
\end{align}
where the summation runs over all
$J:\{1,\ldots,|k|\}\to \{1,\ldots,r\}$ with $\#J^{-1}(i)=k_i$.
(Recall that the variables $t_u$ are identified with
$t^{(j)}_i$'s using (\ref{id}).)
\end{theorem}

\begin{example} Let $r=2$, $k=(2,1)$, and write $t$ for $t^{(1)}$ and $s$ for
$t^{(2)}$. Then
\begin{align}
\Omega_{(2,1)}=&
\Big(\frac{1}{t_1(t_2-t_1)(s-t_2)}+\frac{1}{t_2(t_1-t_2)(s-t_1)}\Big)dt_1\wedge
dt_2\wedge ds \ \otimes \tilde{f}_2\tilde{f}_1^2+\\
\notag &  \Big(\frac{1}{t_1(s-t_1)(t_2-s)}
+\frac{1}{t_2(s-t_2)(t_1-s)}\Big)dt_1 \wedge dt_2 \wedge ds\ \otimes \tilde{f}_1\tilde{f}_2\tilde{f}_1+\\
\notag &  \Big(\frac{1}{s(t_1-s)(t_2-t_1)}+\frac{1}{s(t_2-s)(t_1-t_2)}\Big)
dt_1 \wedge dt_2\wedge ds\  \otimes \tilde{f}_1^2\tilde{f}_2.
\end{align}
\end{example}

Similar rational functions will often appear in this paper. It
will be convenient to encode them with diagrams as follows: $\Omega_{(2,1)}=$
$$
\sym\Big(
\begin{picture}(53,20)
\put(0,0){$*$} \put(3,3){\line(1,0){15}} \put(15,0){$\bullet$}
\put(15,8){$t_1$} \put(18,3){\line(1,0){15}} \put(30,0){$\bullet$}
\put(30,8){$t_2$} \put(33,3){\line(1,0){15}} \put(45,0){$\bullet$}
\put(45,8){$s$}
\end{picture} \Big)dV_k \otimes \tilde{f}_2\tilde{f}_1^2 +
 \sym\Big(
\begin{picture}(53,20)
\put(0,0){$*$} \put(3,3){\line(1,0){15}} \put(15,0){$\bullet$}
\put(15,8){$t_1$} \put(18,3){\line(1,0){15}} \put(30,0){$\bullet$}
\put(30,8){$s$} \put(33,3){\line(1,0){15}} \put(45,0){$\bullet$}
\put(45,8){$t_2$}
\end{picture} \Big)dV_k \otimes \tilde{f}_1\tilde{f}_2\tilde{f}_1+
\sym\Big(
\begin{picture}(53,20)
\put(0,0){$*$} \put(3,3){\line(1,0){15}} \put(15,0){$\bullet$}
\put(15,8){$s$} \put(18,3){\line(1,0){15}} \put(30,0){$\bullet$}
\put(30,8){$t_1$} \put(33,3){\line(1,0){15}} \put(45,0){$\bullet$}
\put(45,8){$t_2$}
\end{picture} \Big)dV_k \otimes \tilde{f}_1^2\tilde{f}_2.
$$
A rooted tree (root denoted by *) with variables associated to its
vertices encodes the product of $1/(a-b)$'s for every edge whose
vertices are decorated by $a$ and $b$, and $b$ is closer to the
root of the tree. The label of the root of the tree is 0, so we do
not write it out. The symmetrizer $\sym$ is meant with respect to
the content of the rational function.

\section{Properties of the differential forms}
In this section we present the two key properties needed in Section \ref{g}.

\subsection{The residue of the canonical differential form.}

For $k=(k_1,\ldots,k_r)$, we denote $(k_1,\ldots,k_{i-1},k_i-1,k_{i+1},\ldots,k_r)$ by
$k-1_i$.

\begin{lemma}\label{l5} Let $k\in\N^r$ and $i\in [1,\ldots, r]$. Then the maps
$$
R\ :\ \A^{G_k} \ \to\ \A^{G_{k-1_i}}\ ,
\qquad
\omega \ \mapsto \ \res_{t^{(i)}_{k_i}=0}\omega\ ,
$$
and
$$
\psi\ :\ U_r[k-1_i]\ \to\ U_r[k]\ ,
\qquad
x \ \mapsto\
(-1)^{k_1+\ldots +  k_i-1} x \tilde{f}_i\ ,
$$
are dual.
\end{lemma}

\begin{proof} Let $\omega\in\A^{G_k}$ and $\tilde{f}_J\in U_r[k-1_i]$. We need to check that
the residue with respect to the flag corresponding to
 $\tilde{f}_J$ of $\res_{t^{(i)}_{k_i}=0} \omega$ is equal to
$(-1)^{k_1+\ldots+k_i-1}$
times the residue with respect to the flag corresponding to $\tilde{f}_J\tilde{f}_i$ of $\omega$.
This follows from the definitions (and the sign conventions).
\end{proof}

\begin{theorem}\label{resth}
\[
\res_{t^{(i)}_{k_i}=0} \Omega_k=(-1)^{k_1+k_2+\ldots+k_i-1}\cdot \Omega_{k-1_i}\cdot(1\otimes \tilde{f}_i).
\]
\end{theorem}

\begin{proof} Let $\{b_u\}$ be a basis of $U_r[k-1_i]$, hence $\Omega_{k-1_i}=\sum b_u^*\otimes b_u$.
Since the map $\psi$ in Lemma \ref{l5} is an embedding,
the images $\psi(b_u)$ can be extended to a basis $\{\psi(b_u), c_v\}$ of $U_r[k]$. Then $\Omega_k=
\sum \psi(b_u)^*\otimes \psi(b_u) + \sum c_v^*\otimes c_v$. We have
$(R\otimes 1)\Omega_k=\sum R(\psi(b_u)^*) \otimes \psi(b_u)+\sum
R(c_v^*)\otimes c_v$, which, according to Lemma \ref{l5}, is $\sum
b_u^*\otimes \psi(b_u)=(1\otimes \psi)\Omega_{k-1_i}$, as required.
\end{proof}

\begin{example}\rm For $r=2$, we write $t$ for $t^{(1)}$ and $s$ for $t^{(2)}$. Then
$\res_{s=0}\Omega_{(1,1)}=$
$$\res_{s=0}\Big(\frac{dt}{t}\wedge\frac{d(s-t)}{s-t}\otimes \tilde{f}_2\tilde{f}_1-
\frac{ds}{s}\wedge\frac{d(t-s)}{t-s}\otimes \tilde{f}_1\tilde{f}_2\Big)=0\otimes \tilde{f}_2\tilde{f}_1
-\frac{dt}{t}\otimes \tilde{f}_1\tilde{f}_2=-\Omega_{(1,0)}\tilde{f}_2.$$
\end{example}

\subsection{The multiplication of differential forms}

Recall that $U_r$ is equipped with a standard Hopf algebra
structure. The co-multiplication $\Delta :U_r\to U_r \otimes U_r$
is defined for degree one elements $x$ as $\Delta(x)=1\otimes
x+x\otimes 1$; e.g. $\Delta(\tilde{f}_1)=1\otimes \tilde{f}_1 +\tilde{f}_1 \otimes 1$.  Then
$\Delta(\tilde{f}_1\tilde{f}_2)=1\otimes \tilde{f}_1\tilde{f}_2+ \tilde{f}_1\otimes \tilde{f}_2 + \tilde{f}_2\otimes \tilde{f}_1
+\tilde{f}_1\tilde{f}_2 \otimes 1$.

The dual $\Delta^*$ of $\Delta$ is therefore a multiplication
on the dual space $U_r^*=\sum_k
\A^{G_k}$. Our goal is to express explicitly this
multiplication of differential forms.

\begin{theorem}\label{starthm} For $k,l\in \N^r$, let
$\omega dV_k\in\A^{G_k}$ and $\eta dV_l\in
\A^{G_l}$ be differential forms.
Then
\begin{equation}\label{deltastar}
\Delta^*(\omega dV_k \otimes \eta dV_l)=(\omega * \eta)\ dV_{k+l},
\end{equation}
(see Section \ref{star}).
\end{theorem}

\begin{proof} We will need the following concept.
Call a triple $(S_1,S_2,J)$ a {\em shuffle} of
$J_1:\{1,\ldots,|k|\}\to \{1,\ldots,r\}$ and $J_2:\{1,\ldots,|l|\}\to
\{1,\ldots,r\}$ if
\begin{itemize}
\item{} $S_1,S_2$ are subsets of $\{1,\ldots,|k+l|\}$,
$\#S_1=|k|$, $\#S_2=|l|$,
\item{}  $\{1,\ldots,|k+l|\}$ is the disjoint union of
$S_1$ and $S_2$,
\item{} $J$ is a map from
$\{1,\ldots,|k+l|\}$ to $\{1,\ldots,r\}$,
\item{} for the increasing bijections $s_1:S_1\to \{1,\ldots,|k|\}$
  and $s_2:S_2\to\{1,\ldots,|l|\}$, we have
$$
J(i)=
\begin{cases}
J_1\circ s_1(i) & i\in S_1\\
J_2\circ s_2(i)& i\in S_2.
\end{cases}$$
\end{itemize}

The collection of $\tilde{f}_{J}$'s form a basis of $U_r$. Let the
dual basis of $U_{r}^*$ be $\{ \tilde{f}_J^*\}$. We only need to
check (\ref{deltastar}) for this dual basis. Hence, let $\omega
dV_k=\tilde{f}_{J_1}^*$, $\eta dV_l=\tilde{f}_{J_2}^*$
 with $\tilde{f}_{J_1}\in U_r[k]$, $\tilde{f}_{J_2}\in U_r[l]$.

The definition of $\Delta$ implies that
$\Delta^*(\tilde{f}_{J_1}^*\otimes \tilde{f}_{J_2}^*)$ is
$\sum \tilde{f}_J^*$, where the summation runs over all shuffles of $J_1$ and $J_2$;
e.g. $\Delta^*(\tilde{f}_1^*\otimes \tilde{f}_2^*)=(\tilde{f}_1\tilde{f}_2)^*+(\tilde{f}_2\tilde{f}_1)^*$,
$\Delta^*(\tilde{f}_1^*\otimes \tilde{f}_1^*)=2(\tilde{f}_1^2)^*$.

On the other hand, the right-hand-side in (\ref{deltastar}) is
also $\sum \tilde{f}_J^*$, with the summation running over the
shuffles of $J_1$ and $J_2$. This can be seen by an iterated
application of the `diagram surgery'
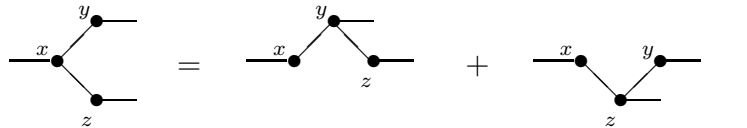
$$\begin{picture}(70,25)(0,3)
\put(0,3){$ $} \put(10,8){\line(1,0){15}} \put(25,5){$\bullet$}
\put(28,8){\line(1,1){15}} \put(28,8){\line(1,-1){15}}
\put(40,20){$\bullet$} \put(40,-10){$\bullet$}
\put(43,23){\line(1,0){15}} \put(43,-7){\line(1,0){15}}
\put(60,20){$ $} \put(60,-10){$ $} \put(20,12){${}_x$}
\put(36,27){${}_y$} \put(37,-15){${}_z$}
\end{picture}=\
\begin{picture}(90,25)(0,3)
\put(0,3){$ $} \put(10,8){\line(1,0){15}} \put(25,5){$\bullet$}
\put(28,8){\line(1,1){15}} \put(42,23){\line(1,-1){15}}
\put(40,20){$\bullet$} \put(55,5){$\bullet$}
\put(43,23){\line(1,0){15}} \put(58,8){\line(1,0){15}}
\put(60,20){$ $} \put(75,5){$ $} \put(20,12){${}_x$}
\put(36,27){${}_y$} \put(53,0){${}_z$}
\end{picture}+\
\begin{picture}(90,25)(0,3)
\put(0,3){$ $} \put(10,8){\line(1,0){15}} \put(25,5){$\bullet$}
\put(43,-7){\line(1,1){15}} \put(28,8){\line(1,-1){15}}
\put(55,5){$\bullet$} \put(40,-10){$\bullet$}
\put(58,8){\line(1,0){15}} \put(43,-7){\line(1,0){15}}
\put(75,5){$ $} \put(60,-10){$ $} \put(20,12){${}_x$}
\put(51,12){${}_y$} \put(37,-15){${}_z$}
\end{picture},$$
\ \vskip 0.5 true cm \noindent justified by the identity
$$\frac{1}{(y-x)(z-x)}=\frac{1}{(y-x)(z-y)}+\frac{1}{(z-x)(y-z)}.$$
c.f. \cite[Lemma~4.4]{mukhin}. For example
$1/(t(s-t))*1/u=1/(t(s-t)u)=$
$$\begin{picture}(40,40)(0,20)
\put(0,20){$*$} \put(3,23){\line(1,1){15}}
\put(15,35){$\bullet$}\put(15,41){$t$} \put(18,38){\line(1,0){15}}
\put(30,35){$\bullet$}\put(35,41){$s$} \put(3,23){\line(1,-1){15}}
\put(15,5){$\bullet$}\put(15,-2){$u$}
\end{picture}=
\begin{picture}(40,40)(0,20)
\put(0,20){$*$} \put(3,23){\line(1,1){15}}
\put(15,35){$\bullet$}\put(15,41){$t$} \put(18,38){\line(1,0){15}}
\put(30,35){$\bullet$}\put(30,41){$s$}
\put(18,38){\line(1,-1){15}}
\put(30,20){$\bullet$}\put(30,14){$u$}
\end{picture}+
\begin{picture}(55,40)(0,20)
\put(0,20){$*$} \put(18,8){\line(1,1){15}}
\put(30,20){$\bullet$}\put(30,26){$t$} \put(33,23){\line(1,0){15}}
\put(45,20){$\bullet$}\put(45,26){$s$} \put(3,23){\line(1,-1){15}}
\put(15,5){$\bullet$}\put(15,-1){$u$}
\end{picture}=
\Big(
\begin{picture}(55,40)(0,20)
\put(0,20){$*$} \put(3,23){\line(1,1){15}}
\put(15,35){$\bullet$}\put(15,41){$t$} \put(18,38){\line(1,0){15}}
\put(30,35){$\bullet$}\put(30,41){$s$}
\put(33,38){\line(1,-1){15}}
\put(45,20){$\bullet$}\put(45,14){$u$}
\end{picture}+
\begin{picture}(55,40)(0,20)
\put(0,20){$*$} \put(3,23){\line(1,1){15}}
\put(15,35){$\bullet$}\put(15,41){$t$} \put(33,23){\line(1,0){15}}
\put(45,20){$\bullet$}\put(45,26){$s$}
\put(18,38){\line(1,-1){15}}
\put(30,20){$\bullet$}\put(30,14){$u$}
\end{picture}\Big)+
\begin{picture}(55,40)(0,20)
\put(0,20){$*$} \put(18,8){\line(1,1){15}}
\put(30,20){$\bullet$}\put(30,26){$t$} \put(33,23){\line(1,0){15}}
\put(45,20){$\bullet$}\put(45,26){$s$} \put(3,23){\line(1,-1){15}}
\put(15,5){$\bullet$}\put(15,-1){$u$}
\end{picture}=
$$
\ \vskip .5 true cm \noindent
$=1/(t(s-t)(u-s))+1/(t(u-t)(s-u))+1/(u(t-u)(s-t))$.
\end{proof}

\section{The canonical differential form for the simple Lie algebras $A,B,C,D$}
\label{g}

The following projections of the canonical differential form are used
in integral solutions of the KZ equations and in the Bethe ansatz construction. Let $\g$
be a simple Lie algebra of rank $r$, with Cartan decomposition
$\g=\n\oplus \h \oplus \np$. The
universal enveloping algebra $U(\n)$ of $\n$ is generated by $r$
elements $f_1,\ldots, f_r$ (the standard Chevalley generators)
subject to the Serre relations; i.e. there is the quotient map $q:
U_r\to U(\n)$ sending $\tilde{f_i}$ to $f_i$ for any $i$. We say
that an element $x \in U(\n)$ has content $k \in \N^r$ if $x \in
q(U_r[k])$.


\begin{definition}
The canonical differential form $\Omega^{\g}_{k}$ of a simple
Lie algebra $\g$ is defined as the image of $\Omega_{k}$ under the map
$$\id\otimes\ q:
\A^{G_k}\otimes U_r[k] \to
\A^{G_k}\otimes U(\n)[k].$$
\end{definition}

The Lie algebra $\n$ is a direct sum of $1$-dimensional weight
spaces $\mathfrak{n}_\beta$ labelled by the positive roots
$\beta:\,\n=\oplus_{\beta} \mathfrak{n}_\beta$.  Let $F_{\beta}$ be
a choice of generator in $\mathfrak{n}_\beta$.  We fix a linear
ordering of the positive roots: $\beta_1,\dots,\beta_m$. According
to the Poincar\'e-Birkhoff-Witt theorem, a $\CC$-basis of the
algebra $U(\n)$ is given by the collection of elements
$F_{\beta_1}^{p_1}F_{\beta_2}^{p_2}\ldots F_{\beta_m}^{p_m}$, where
$m=\dim \n$ and $p=(p_1,\ldots,p_m)\in\N^m$.

\begin{example} \rm
For $sl_3$, the Lie algebra of all $3\times 3$ matrices with trace equal to zero, the
positive roots $\alpha_1$, $\alpha_2$, $\alpha_1+\alpha_2$ correspond
to the matrix entries at the positions $(2,1)$, $(3,2)$, and $(3,1)$,
and in turn, to the basis $F_{\alpha_1}=f_1, F_{\alpha_2}=f_2,
F_{\alpha_1+\alpha_2}=[f_2,f_1]$ of $\n$.  A PBW basis of $U(\n)$ is
$f_1^{p_1}[f_2,f_1]^{p_2}f_2^{p_3}$ with $p=(p_1,p_2,p_3)\in \N^3$.
\end{example}

After fixing the linear ordering $\beta_1<\beta_2<\ldots<\beta_m$
of the positive roots of $\g$, the canonical differential form of
$\g$ can be written in the form of
\begin{equation}
\label{og}\Omega_k^{\g}\ =\ \sum_{p} \ \omega_p \ dV_k  \otimes
  F_{\beta_1}^{p_1}F_{\beta_2}^{p_2}\ldots
F_{\beta_m}^{p_m}\  .
\end{equation}
Here the summation is over $p$ such that the content of
$F_{\beta_1}^{p_1}F_{\beta_2}^{p_2}\ldots
F_{\beta_m}^{p_m}$ is $k$; \
 $\omega_p dV_k$ is a differential form in $\A^{G_k}$,
and $\omega_p$ is a rational function.

\begin{theorem}{\rm (Product formula.)} For $l = 1, \dots , m$, let the content of
  $F_{\beta_l}$ be $k^{(l)}$. Then
there exist rational functions $\eta_{\beta_l}$ in the variables
$(t^{(i)}_j)_{j=1,\ldots,k^{(l)}_i}$, symmetric under
$G_{k^{(l)}}$, such that
$$\omega_p=\frac{1}{\prod_l p_l!}\cdot\overbrace{\eta_{\beta_1}*\ldots*\eta_{\beta_1}}^{p_1}*
\overbrace{\eta_{\beta_2}*\ldots*\eta_{\beta_2}}^{p_2}* \ldots*
\overbrace{\eta_{\beta_m}*\ldots*\eta_{\beta_m}}^{p_m}$$
\end{theorem}

\begin{proof} Denote $F^p=F_{\beta_1}^{p_1}\ldots F_{\beta_m}^{p_m}$. The co-multiplication
$\Delta$
can be expressed in the PBW basis as:
$$
\Delta(F^p)=\sum_{p'+p''=p}\ \
\prod_{i=1}^m \frac{p_i!}{p'_i!p''_i!} \cdot F^{p'}\otimes F^{p''}\ .
$$
Then for the dual multiplication we have
$$
\Delta^*(\omega_{p'}dV_k\otimes \omega_{p''}dV_l) =
\Delta^*(F^{p'*}\otimes F^{p''*})
=
\left(\prod_i \frac{(p'_i+p''_i)!}{p'_i!p''_i!}\right)\cdot F^{(p'+p'')*}.
$$
Using Theorem \ref{starthm}, we obtain
$$\omega_{p'}* \omega_{p''}=\prod_i \frac{(p'_i+p''_i)!}{p'_i!p''_i!}\cdot\omega_{p'+p''},$$
from which the result follows (put $\eta_{\beta_l}=\omega_{1_l}$).
\end{proof}

\begin{example} \rm For $\g=sl_3$ and the ordering $\alpha_1<\alpha_1+\alpha_2<\alpha_2$,
we have $\omega_{\alpha_1}=1/t_1$, $\omega_{\alpha_1+\alpha_2}=1/(t_1(s_1-t_1))$ and
$\omega_{\alpha_2}=1/s_1$. (Again, we write $t$ for $t^{(1)}$ and $s$ for $t^{(2)}$.)
Then the differential form corresponding to $f_1^2[f_2,f_1]$ is
$$\omega_{(2,1,0)}dt_1\wedge dt_2\wedge dt_3\wedge ds=
\frac{1}{2!1!0!} \cdot \frac{1}{t_1}*\frac{1}{t_1}*\frac{1}{t_1(s_1-t_1)}=$$
$$=\frac{1}{2} \sym_{(3,1)}\Big(\frac{1}{t_1t_2t_3(s-t_3)}\Big)dV_{(3,1)}=
\frac{3s^2-2s(t_1+t_2+t_3)+(t_1t_2+t_1t_3+t_2t_3)}{t_1t_2t_3(s-t_1)(s-t_2)(s-t_3)}dV_{(3,1)}.$$
\end{example}

\bigskip

This means that $\Omega_k^{\g}$ is determined once we know its
`atoms', i.e. the $\eta_{\beta}$'s for the positive roots $\beta$.
In the remainder of this section, we will compute them for the infinite series $A$, $B$, $C$, $D$
of simple Lie algebras. For each of these, we will list (1) the positive
roots,
(2) the simple roots, and (3) the expression of the positive
roots in terms of the simple roots.
Then we choose (4) a linear ordering of the positive roots
and fix (5) the elements
$F_{\beta}$'s (choice of a constant). Then we
describe the elements $\eta_{\beta}$'s with the choices (4), (5).

Let $\epsilon_i=(0,\dots,0,1,0,\dots,0)\in\CC^{r}$ ($1$ occurs at
the $i$th position from the left). For a multi-index
$J=(J(1),J(2),\ldots,J(n))$, let
$[f_J]=[f_{J(1)},[f_{J(2)},[\ldots,[f_{J(n-1)},f_{J(n)}]\ldots]]]$.

\subsection{The simple Lie algebra $A_{r-1}$} \label{pbwa}
\begin{enumerate}
\item{} The positive roots are $\epsilon_i-\epsilon_j$ for $1\leq i<j\leq r$.
\item{} The simple roots are $\alpha_i=\epsilon_i-\epsilon_{i+1}$ for $1\leq i< r$.
\item{} $\epsilon_i-\epsilon_j=\sum_{u=i}^{j-1}\alpha_u$.
\item{} Let $\epsilon_i-\epsilon_j<\epsilon_{i'}-\epsilon_{j'}$ if either $i<i'$, or $i=i'$ and $j<j'$.
\item{} For a positive root $\beta=\epsilon_i-\epsilon_j$, let
$F_{\beta}=[f_{(j-1,j-2,\ldots,i)}]$.
\end{enumerate}
\begin{theorem} \label{ar} For
$\beta=\epsilon_i-\epsilon_j=\alpha_i+\alpha_{i+1}+\ldots+\alpha_{j-1}$,
we have
$$\eta_{\beta}=\frac{1}{t_1^{(i)}(t_1^{(i+1)}-t_1^{(i)})\cdot \ldots\cdot (t_1^{(j-1)}-t_1^{(j-2)})}.$$
\end{theorem}
The result can be visualized by the following string-diagram (the
labels of the vertices in the diagram indicate the superscripts of
the corresponding $t_1$'s).
$$\eta_{\alpha_i+\ldots+\alpha_j}=
\begin{picture}(85,20)
\put(0,0){$*$} \put(0,3){\line(1,0){15}} \put(15,0){$\bullet$}
\put(15,10){${}_i$} \put(15,3){\line(1,0){15}} \put(30,0){$\bullet$}
\put(28,10){${}_{i+1}$} \put(30,3){\line(1,0){15}} \put(50,0){$\ldots$}
\put(65,3){\line(1,0){15}} \put(80,0){$\bullet$} \put(80,10){${}_{j-1}$}
\end{picture}$$

\subsection{The simple Lie algebra $B_r$}
\begin{enumerate}
\item{} The positive roots are $\epsilon_i$ ($1\leq i\leq r$) and
$\epsilon_i-\epsilon_j$, $\epsilon_i+\epsilon_j$ ($1\leq i<j\leq
r$). \item{} The simple roots are
$\alpha_i=\epsilon_i-\epsilon_{i+1}$ for $i=1,\ldots,r-1$ and
$\alpha_r=\epsilon_r$ (the `short' root). \item{}
$$\epsilon_i=\sum_{u=i}^r \alpha_u \ (1\leq i \leq r), \quad
\epsilon_i-\epsilon_j=\sum_{u=i}^{j-1} \alpha_u, \quad
\epsilon_i+\epsilon_j=\sum_{u=i}^{j-1} \alpha_u
+2\sum_{u=j}^{r}\alpha_u\ (1\leq i<j\leq r).$$ \item{} Let $\beta$
be one of $\epsilon_i$, $\epsilon_i-\epsilon_j$, or
$\epsilon_i+\epsilon_j$,  and let $\beta'$ be one of
$\epsilon_{i'}$, $\epsilon_{i'}-\epsilon_{j'}$, or
$\epsilon_{i'}+\epsilon_{j'}$ Then we set $\beta<\beta'$ if
$i>i'$. If $i<j<j'$ then we also set
$\epsilon_i+\epsilon_j<\epsilon_i+\epsilon_{j'}<\epsilon_i<\epsilon_i-\epsilon_{j'}<\epsilon_i-\epsilon_j$.
\item{}
$F_{\alpha_i+\ldots+\alpha_{j-1}}=[f_{(j-1,j-2,\ldots,i)}]$ and
$F_{\epsilon_i+\epsilon_j}=[
[f_{(i,i+1,\ldots,r)}],[f_{(r,r-1,\ldots,j)}]]$.
\end{enumerate}

\noindent The vector $(\overbrace{0,\ldots,
  0}^{u_0},\overbrace{1,\ldots,1}^{u_1}\dots)$ will be abbreviated by $(0^{u_0}1^{u_1}\ldots)$.

\begin{theorem}
We have
$$\eta_{\alpha_i+\ldots+\alpha_{j-1}}=\frac{1}{t_1^{(j-1)}(t_1^{(j-2)}-t_1^{(j-1)})\cdot
\ldots\cdot (t_1^{(i)}-t_1^{(i+1)})} \qquad\qquad(\text{for both roots}\ \epsilon_i\ \text{and}\ \epsilon_i-\epsilon_j),$$
\begin{align}\notag
\eta_{\epsilon_i+\epsilon_j}=\frac{1}{2}\sym_{(0^{i-1}1^{j-i}2^{r-j+1})}\Big(&
\frac{t_1^{(r-1)}-t_2^{(r-1)}}{
t_1^{(j)}(t_1^{(r)}-t_1^{(r-1)})(t_2^{(r)}-t_1^{(r-1)})(t_2^{(r-1)}-t_1^{(r)})(t_2^{(r-1)}-t_2^{(r)})}\\ \notag
& \cdot\frac{1}{\prod_{k=j+1}^{r-1}(t_1^{(k)}-t_1^{(k-1)}) \prod_{k=i}^{r-2}(t_2^{(k)}-t_2^{(k+1)}) }\Big).
\end{align}
\end{theorem}\label{br}
The structure of these functions is better understood via the following pictures.
$$\eta_{\alpha_i+\ldots+\alpha_{j-1}}=
\begin{picture}(85,20)
\put(0,0){$*$} \put(3,3){\line(1,0){15}} \put(15,0){$\bullet$}
\put(10,10){${}_{j-1}$} \put(18,3){\line(1,0){15}} \put(30,0){$\bullet$}
\put(28,10){${}_{j-2}$} \put(33,3){\line(1,0){15}} \put(50,0){$\ldots$}
\put(68,3){\line(1,0){15}} \put(80,0){$\bullet$} \put(82,10){${}_i$}
\end{picture}$$

\medskip

$$\begin{picture}(190,40)(-70,-30)
\put(-170,0){$\eta_{\sum_{i\leq k<j} \alpha_k +2\sum_{j\leq k \leq r}\alpha_k}=\frac{1}{2}\sym$}
\put(0,0){$*$}\put(3,3){\line(1,0){15}}
\put(15,0){$\bullet$} \put(15,9){${}_j$}
\put(18,3){\line(1,0){15}}
\put(30,0){$\bullet$} \put(28,9){${}_{j+1}$}
\put(33,3){\line(1,0){15}}
\put(50,0){$\ldots$} \put(68,3){\line(1,0){15}}
\put(80,0){$\bullet$} \put(74,10){${}_{r-1}$}
\put(83,2){\line(1,0){30}}
\put(83,4){\line(1,0){30}}
\put(83,3){\line(1,1){15}}
\put(83,3){\line(1,-1){15}}
\put(95,15){$\bullet$} \put(95,25){${}_r$}
\put(95,-15){$\bullet$} \put(95,-20){${}_r$}
\put(98,19){\line(1,-1){15}}
\put(98,-13){\line(1,1){15}}
\put(110,0){$\bullet$} \put(110,-6){${}_{r-1}$}
\put(113,3){\line(1,0){15}}
\put(125,0){$\bullet$}\put(125,8){${}_{r-2}$}
\put(128,3){\line(1,0){15}}
\put(145,0){$\ldots$}
\put(168,3){\line(1,0){15}}
\put(180,0){$\bullet$}\put(182,8){${}_i$}
\end{picture}$$
In the second picture, the double edge means that the
corresponding difference is in the numerator. The labels mean
superscripts of variables $t_1$, and when a superscript $i$ is
used twice in a diagram, they mean $t_1^{(i)}$ and $t_2^{(i)}$.

\subsection{The simple Lie algebra $C_r$}
\begin{enumerate}
\item{} The positive roots
are $\epsilon_i-\epsilon_j$,
$\epsilon_i+\epsilon_j$ for $1\leq i<j\leq r$ and $2\epsilon_i$ for $1\leq i\leq r$.
\item{} The simple roots are $\alpha_i=\epsilon_i-\epsilon_{i+1}$ ($1\leq i<r$) and
 $\alpha_r=2\epsilon_r$
(the `long' root).
\item{}
$$\epsilon_i-\epsilon_j=\sum_{u=i}^{j-1}\alpha_u, \qquad
\epsilon_i+\epsilon_j=\sum_{u=i}^{j-1}\alpha_u+2\sum_{u=j}^{r-1}\alpha_u+\alpha_r,
\qquad 2\epsilon_i=2\sum_{u=i}^{r-1}\alpha_u+\alpha_r.$$ \item{}
Let $\beta$ be one of $\epsilon_i-\epsilon_j$,
$\epsilon_i+\epsilon_j$, or $2\epsilon_i$, and let $\beta'$ be one
of $\epsilon_{i'}-\epsilon_{j'}$, $\epsilon_{i'}+\epsilon_{j'}$,
or $2\epsilon_{i'}$. Then we set $\beta<\beta'$ if $i>i'$. For
$i<j<j'$ we also set
$\epsilon_i+\epsilon_j<\epsilon_i+\epsilon_{j'}<2\epsilon_i<\epsilon_i-\epsilon_{j'}<\epsilon_i-\epsilon_j$.
\item{} $F_{\epsilon_i-\epsilon_j}=[f_{(i,i+1,\ldots,j-1)}]$,
$F_{\epsilon_i+\epsilon_j}=[f_{i,i+1,\ldots,r-1,r,r-1,\ldots,j}]$,
$F_{2\epsilon_i}=[[f_{(i,i+1,\ldots,r-1)}],[f_{(i,i+1,\ldots,r)}]]$.
\end{enumerate}
\begin{theorem}\label{cr} We have
\begin{align}
\eta_{\epsilon_i-\epsilon_j}=&\frac{1}{t_1^{(j-1)}\prod_{u=i}^{j-2}(t_1^{(u)}-t_1^{(u+1)})},  \notag\\
\eta_{\epsilon_i+\epsilon_j}=&\sym_{(0^{i-1}1^{j-i}2^{r-j}1)}\Big(\ \frac{1}{t_1^{(j-1)}\prod_{u=j+1}^{r}(t_1^{(u)}-t_1^{(u-1)})
\prod_{u=i}^{r-1}(t_2^{(u)}-t_2^{(u+1)})}\ \!\!\Big),\notag
\\
{}\!\!\!
\eta_{2\epsilon_i}=& \sym_{(0^{i-1}2^{r-i}1)}\Big( \ \frac{1}{t_1^{(r)}(t_1^{(r-1)}-t_1^{(r)})(t_2^{(r-1)}-t_1^{(r)})
\prod_{u=i}^{r-2}(t_1^{(u)}-t_1^{(u+1)})\prod_{u=i}^{r-2}(t_2^{(u)}-t_2^{(u+1)})} \ \Big).
\notag\end{align}
\end{theorem}

The result can be visualized by the following diagrams (labels mean upper indices).
\begin{align}
\eta_{\epsilon_i-\epsilon_j}= &
\begin{picture}(85,20)
\put(0,0){$*$} \put(3,3){\line(1,0){15}} \put(15,0){$\bullet$}
\put(10,10){${}_{j-1}$} \put(18,3){\line(1,0){15}} \put(30,0){$\bullet$}
\put(27,10){${}_{j-2}$} \put(33,3){\line(1,0){15}} \put(52,2){$\ldots$}
\put(68,3){\line(1,0){15}} \put(80,0){$\bullet$} \put(80,10){${}_i$}
\end{picture} \notag\\
\eta_{\epsilon_i+\epsilon_j}=& \sym
\begin{picture}(185,20)
\put(0,0){$*$}\put(3,3){\line(1,0){15}}
\put(15,0){$\bullet$} \put(15,10){${}_j$}\put(18,3){\line(1,0){15}}
\put(30,0){$\bullet$} \put(28,10){${}_{j+1}$} \put(33,3){\line(1,0){15}}
\put(52,2){$\ldots$} \put(68,3){\line(1,0){15}}
\put(80,0){$\bullet$} \put(74,10){${}_{r-1}$} \put(83,3){\line(1,0){15}}
\put(95,0){$\bullet$} \put(95,10){${}_r$} \put(95,3){\line(1,0){15}}
\put(110,0){$\bullet$} \put(105,10){${}_{r-1}$}
\put(113,3){\line(1,0){15}}
\put(125,0){$\bullet$}\put(125,10){${}_{r-2}$}
\put(128,3){\line(1,0){15}}
\put(147,2){$\ldots$}
\put(168,3){\line(1,0){15}}
\put(180,0){$\bullet$}\put(182,10){${}_i$}
\end{picture}\notag\\
\eta_{2\epsilon_i}=&\sym
\begin{picture}(150,35)(0,15)
\put(0,15){$*$} \put(3,18){\line(1,0){15}}
\put(15,15){$\bullet$} \put(18,18){\line(1,1){15}} \put(18,18){\line(1,-1){15}}
\put(23,17){${}_r$}
\put(30,30){$\bullet$}  \put(33,33){\line(1,0){15}} \put(24,39){${}_{r-1}$}
\put(30,0){$\bullet$} \put(33,3){\line(1,0){15}} \put(24,-6){${}_{r-1}$}
\put(45,30){$\bullet$} \put(48,33){\line(1,0){15}} \put(45,39){${}_{r-2}$}
\put(45,0){$\bullet$} \put(48,3){\line(1,0){15}} \put(45,-6){${}_{r-2}$}
\put(70,32){$\ldots$} \put(98,33){\line(1,0){15}}
\put(70,2){$\ldots$} \put(98,3){\line(1,0){15}}
\put(110,30){$\bullet$} \put(110,39){${}_i$}
\put(110,0){$\bullet$} \put(110,-6){${}_i$}
\end{picture}
\notag
\end{align}

\bigskip
\bigskip

\subsection{The simple Lie algebra $D_r$}
\begin{enumerate}
\item{} The positive roots
are $\epsilon_j-\epsilon_i$ and $\epsilon_j+\epsilon_i$ for $1\leq i<j\leq r$.
\item{}The simple roots are $\alpha_1=\epsilon_1+\epsilon_2$ and $\alpha_i=\epsilon_i-\epsilon_{i-1}$ ($1<i\leq r$).
\item{}
$$\epsilon_j-\epsilon_i=\sum_{u=i+1}^{j}\alpha_u, \quad
\epsilon_i+\epsilon_j=\alpha_1+\alpha_2+2\sum_{u=3}^{i}\alpha_u+\sum_{u=i+1}^{j}\alpha_u\
(1<i<j\leq r), \quad
\epsilon_1+\epsilon_j=\alpha_1+\sum_{u=3}^{j}\alpha_u.$$ \item{}
Let $\beta$ be one of $\epsilon_j-\epsilon_i$ or
$\epsilon_j+\epsilon_i$, and let $\beta'$ be one of
$\epsilon_{j'}-\epsilon_{i'}$ or $\epsilon_{j'}+\epsilon_{i'}$.
Then we set $\beta<\beta'$ if $j<j'$. For $i'<i<j$ we also set
$\epsilon_j+\epsilon_i<\epsilon_j+\epsilon_{i'}<\epsilon_j-\epsilon_{i'}<\epsilon_j-\epsilon_i$.
\item{} $F_{\epsilon_j-\epsilon_i}=[f_{(j,j-1,\ldots,i+1)}]$,
$F_{\epsilon_1+\epsilon_j}=[f_{(j,j-1,\ldots,3,1)}]$,
$F_{\epsilon_i+\epsilon_j}=[f_{(j,j-1,\ldots,2,1,3,4,\ldots,i-1,i)}]$.
\end{enumerate}
\begin{theorem}\label{dr}We have
\begin{align}
\eta_{\epsilon_j-\epsilon_i}=&\frac{1}{t_1^{(i+1)}\prod_{u=i+2}^{j}(t_1^{(u)}-t_1^{(u-1)})},\notag\\
\eta_{\epsilon_j+\epsilon_1}=&\frac{1}{t_1^{(1)}\prod_{u=3}^{j}(t_1^{(u)}-t_1^{(u-1)})}.\notag
\end{align}
For $1<i<j$, $\eta_{\epsilon_j+\epsilon_i}=\sym_{(1^22^{i-2}1^{j-i}0^{r-j})}$
$$
\frac{t_1^{(3)}-t_2^{(3)}}{
t_2^{(i)}(t_1^{(1)}-t_2^{(3)})(t_1^{(2)}-t_2^{(3)})(t_1^{(3)}-t_1^{(1)})(t_1^{(3)}-t_1^{(2)})
\prod_{u=3}^{i+1}(t_2^{(u)}-t_2^{(u+1)}) \prod_{u=4}^{j}(t_1^{(u)}-t_1^{(u-1)}) }.$$
\end{theorem}

The result can be visualized by the following diagrams (labels mean upper indices).
\begin{align}
\eta_{\epsilon_j-\epsilon_i}=&
\begin{picture}(85,20)
\put(0,0){$*$} \put(3,3){\line(1,0){15}} \put(15,0){$\bullet$}
\put(12,10){${}_{i+1}$} \put(18,3){\line(1,0){15}}
\put(30,0){$\bullet$} \put(28,10){${}_{i+2}$}
\put(33,3){\line(1,0){15}} \put(52,2){$\ldots$}
\put(68,3){\line(1,0){15}} \put(80,0){$\bullet$}
\put(80,10){${}_j$}
\end{picture}\notag\\
\eta_{\epsilon_j+\epsilon_1}=&
\begin{picture}(85,20)
\put(0,0){$*$} \put(3,3){\line(1,0){15}} \put(15,0){$\bullet$}
\put(15,10){${}_{1}$} \put(18,3){\line(1,0){15}} \put(30,0){$\bullet$}
\put(30,10){${}_{3}$} \put(33,3){\line(1,0){15}} \put(52,2){$\ldots$}
\put(68,3){\line(1,0){15}} \put(80,0){$\bullet$} \put(80,10){${}_j$}
\end{picture}\notag\\
\eta_{\epsilon_j+\epsilon_i}=&\sym_k \begin{picture}(190,20)
\put(0,0){$*$}\put(3,3){\line(1,0){15}}
\put(15,0){$\bullet$} \put(15,9){${}_i$}
\put(18,3){\line(1,0){15}}
\put(30,0){$\bullet$} \put(28,9){${}_{i-1}$}
\put(33,3){\line(1,0){15}}
\put(50,0){$\ldots$} \put(68,3){\line(1,0){15}}
\put(80,0){$\bullet$} \put(76,10){${}_{3}$}
\put(83,2){\line(1,0){30}}
\put(83,4){\line(1,0){30}}
\put(83,3){\line(1,1){15}}
\put(83,3){\line(1,-1){15}}
\put(95,15){$\bullet$} \put(95,25){${}_1$}
\put(95,-15){$\bullet$} \put(95,-20){${}_2$}
\put(98,19){\line(1,-1){15}}
\put(98,-13){\line(1,1){15}}
\put(110,0){$\bullet$} \put(112,-6){${}_{3}$}
\put(113,3){\line(1,0){15}}
\put(125,0){$\bullet$}\put(125,8){${}_{4}$}
\put(128,3){\line(1,0){15}}
\put(145,0){$\ldots$}
\put(168,3){\line(1,0){15}}
\put(180,0){$\bullet$}\put(182,9){${}_j$}
\end{picture}\notag
\end{align}

\bigskip\medskip

\subsection{The proofs of the theorems \ref{ar}-\ref{dr}}
Let $\g$ be one of the simple Lie algebras $A_r$, $B_r$, $C_r$,
$D_r$, and let $\beta$ be one of its positive roots. In one of the
Theorems \ref{ar}--\ref{dr} (the one referring to $\g$), we state a
formula for $\eta_\beta$; let us denote the function on
the right-hand-side of that formula by $\overline{\eta}_\beta$. In
this section we will prove that
$\eta_\beta=\overline{\eta}_\beta$, by proving Theorems
\ref{lem1} and \ref{lem2} below.

\begin{lemma} \label{l1}
Under the correspondence between $U_r[k]$ and $Fl^{G_k}$ of
Section \ref{flags}, we have
\begin{align} \notag [\tilde f_i,\tilde f_j] &
\leftrightarrow  \pm  [\CC^{|k|}\supset
(t_1^{(i)}=t_1^{(j)})\supset (t_1^{(i)}=t_1^{(j)}=0)] \\
\notag [\tilde f_i,[\tilde f_i,\tilde f_j]] & \leftrightarrow
\pm\asym_{(2,1)}\Big( [\CC^{|k|}\supset (t^{(i)}_1=t^{(j)}_1)
\supset (t^{(i)}_1=t^{(i)}_2=t^{(j)}_1)
\supset (t^{(i)}_1=t^{(i)}_2=t^{(j)}_1=0)]\ \Big) \\
\notag [\tilde f_i,[\tilde f_i,[\tilde f_i,\tilde f_j]]] &
\leftrightarrow {}
\end{align}
$$\ \qquad\qquad \pm\asym_{(3,1)}\ \Big( [\CC^{|k|}\supset (t^{(i)}_1=t^{(j)}_1)
\supset (t^{(i)}_1=t^{(j)}_1=t^{(i)}_2) \supset
(t^{(i)}_1=t^{(j)}_1=t^{(i)}_2=t^{(i)}_3)$$ $$\qquad\qquad \supset
(t^{(i)}_1=t^{(j)}_1=t^{(i)}_2=t^{(i)}_3=0)] \ \Big).$$
\end{lemma}

\begin{proof}
For $i<j$ let $L=(t^{(i)}_1=t^{(j)}_1=0)$. Then we have
$$
\tilde f_i\tilde f_j -\tilde f_j\tilde f_i
 \leftrightarrow -[\CC^{|k|}\supset (t^{(j)}_1=0)\supset L] -[\CC^{|k|}\supset
(t^{(i)}_1=0)\supset L]= [\CC^{|k|}\supset
(t^{(j)}_1=t^{(i)}_1)\supset L],$$ which proves the first
statement. The others follow from similar calculations.
\end{proof}

Let the content of $F_{\beta}$ be $k$.

\begin{theorem}\label{lem1}
Let $F\in Fl^{G_k}$ be a linear combination of flags corresponding
to an element $\sum c_J \tilde f_J$ in $U_r[k]$ of content $k$. If $\sum_j
c_J\tilde f_j$ belongs to the ideal generated by the Serre relations,
then $\res_F \overline{\eta}_\beta dV_k=0$.
\end{theorem}

\begin{proof}
We show the result for $\g$ of type $A$. In this case there are two
kinds of Serre relations: $[f_i,f_j]=0$ if
$(\alpha_i,\alpha_j)=0$, and $[f_i,[f_i,f_j]]=0$ if
$(\alpha_i,\alpha_j)=-1$. We will consider the linear combination
of flags corresponding to multiples of $[\tilde f_i,\tilde f_j]$
and $[\tilde f_i,[\tilde f_i,\tilde f_j]]$. By Lemma~\ref{l1}, any
multiple of $[\tilde f_i,\tilde f_j]$ corresponds to a linear
combination $F$ of flags of the form $\pm
[\CC^{|k|}\supset\ldots\supset L_{u+2}\supset L_{u+1} \supset L_u
\supset \ldots \supset 0]$ with
$$
[L_{u+2}/L_u\supset L_{u+1}/L_u \supset L_u/L_u] \ \simeq\
[\CC^2\supset (t^{(i)}_{v_1}=t^{(j)}_{v_2}) \supset 0],
$$
where the coordinates on $\CC^2$ are $t^{(i)}_{v_1}$ and
$t^{(j)}_{v_2}$.

The rational function $\overline{\eta}_\beta$ does not have a
factor of type $t^{(i)}-t^{(j)}$ in the denominator for
$(\alpha_i,\alpha_j)=0$. Thus $\res_{F}\overline{\eta}_\beta
dV_k=0$.

By Lemma \ref{l1},  any multiple of $[\tilde f_i[\tilde f_i,\tilde
f_j]]$ corresponds to a linear combination $F$ of flags, with each
term of the form $\pm [\CC^{|k|}\supset\ldots\supset
L_{u+3}\supset L_{u+2}\supset L_{u+1} \supset L_u \supset \ldots
\supset 0]$ with
$$
[L_{u+3}/L_u\supset L_{u+2}/L_u\supset L_{u+1}/L_u \supset
L_u/L_u] \ \simeq\ [\CC^3\supset (t_{v_1}^{(i)}-t_{v_2}^{(j)})
\supset (t^{(i)}_{v_1}=t^{(i)}_{v_3}=t_{v_2}^{(j)}) \supset 0],$$
where the coordinates on $\CC^3$ are $t^{(i)}_{v_1}$,
$t^{(i)}_{v_3}$, and $t_{v_2}^{(j)}$ with $v_1\not=v_3$. Since all
$k_i$ are 0 or 1 for the positive root $\beta$, this type of flags
cannot occur in $Fl^{G_k}$.

The proof for the types $B, C, D$ are analogous.
\end{proof}

\begin{remark} \rm
By Theorem \ref{lem1},
the differential forms $\Omega^{\g}$ (and $\Omega^{V}$, see the
Introduction and Section \ref{repval} below) do not have poles at
the $t^{(i)}=t^{(j)}$ type hyperplanes if the corresponding simple
roots are orthogonal. Hence the poles of $\Omega^{\g}$ coincide
with the singularities of the master function, see the
Introduction.
\end{remark}

\medskip

Now let $F^p=F_{\beta_1}^{p_1}\ldots F_{\beta_m}^{p_m}$ be another
element of content $k$ in the Poincare-Birkhoff-Witt basis,
different from $F_{\beta}$. Let us choose preimages of $F_{\beta}$
and $F^p$ under the projection $U_r[k]\to U(\n)$, and let
$Flag_\beta$ and $Flag^p\in Fl^{G_k}$ be the corresponding linear
combinations of flags.

\begin{theorem} \label{lem2} \
\begin{itemize}
\item{} The residue of the differential form
$\overline{\eta}_\beta dV_k$ with respect to $Flag^p$ is 0.
\item{} The residue of the differential form
$\overline{\eta}_\beta dV_k$ with respect to $Flag_\beta$ is 1.
\end{itemize}
\end{theorem}

\begin{lemma} \label{resres}
For $i=1,\ldots,r$, $j=1,\ldots, k_i$, we have
$$
\res_{t^{(i)}_{j}=0}\, \eta_\beta\,dV_k\ =\
\res_{t^{(i)}_{j}=0}\,\overline{\eta}_\beta\,dV_k\ .
$$
\end{lemma}

\begin{proof}
It is enough to consider $j=k_i$. Then the left hand side can be
calculated from Theorem~\ref{resth}, and the right hand side is
given explicitly. For types $A, B, D$, we obtain
$$
\res_{t_1^{(i)}=0} \eta_\beta dV_k=\res_{t_1^{(i)}=0}
\overline{\eta}_\beta dV_k\ =\
\begin{cases} \eta_{\beta-\alpha_i}dV_{k-1_i}
 & \hbox{if}\ \beta-\alpha_i\ \hbox{is a positive root},\ \beta-\alpha_i>\beta
\\
0 & \hbox{otherwise.}\end{cases}$$ For type $C$ we obtain $
\res_{t_1^{(i)}=0} \eta_\beta dV_k=\res_{t_1^{(i)}=0}
\overline{\eta}_\beta dV_k =$
$$=\begin{cases}
\eta_{\beta-\alpha_i}dV_{k-1_i}
 & \hbox{if}\ \beta-\alpha_i\ \hbox{is a positive root},\ \beta-\alpha_i>\beta
\\
\eta_{\frac{\beta-\alpha_i}{2}}*\eta_{\frac{\beta-\alpha_i}{2}}\ d
V_{k-1_i} & \hbox{if}\ \frac{\beta-\alpha_i}{2}\ \hbox{is a
positive root}, \frac{\beta-\alpha_i}{2}>\beta
\\
0 & \hbox{otherwise.}\end{cases}$$
\end{proof}

\noindent{\sl Proof of Theorem~\ref{lem2}.} The second statement
follows from the explicit forms for $\overline{\eta}_\beta$.

Let $\overline{\Omega}^{\g}_k$ be the form obtained from
$\Omega^{\g}_k$ by replacing the term $\eta_\beta\otimes F_\beta$
by $\overline{\eta}_\beta \otimes F_\beta$. Then $\res_{Flag^p}
\overline{\Omega}^{\g}_{k}=\res_{Flag^p} {\Omega}^{\g}_{k}$
by Lemma \ref{resres}.  We have $\res_{Flag^p} {\Omega}^{\g}_{k}=1\otimes F^p$
by the fact $(U_r[k])^*=\A^{G_k}$. Therefore we have $\res_{Flag^p}
\overline{\eta}_\beta=0$, as required. \qed

\section{ Appendix: Representation-valued canonical differential form}
\label{repval}
For the convenience of the reader, we give formulas from \cite{sv1}, \cite{sv} for
the $V$-valued differential form $\Omega^V$ which appears in the hypergeometric
solutions to the KZ equations and the Bethe ansatz method.

Consider the formula for the canonical differential form from Theorem
\ref{canform}, without putting $t_0=0$. Instead, put $t_0=z$ and
denote this form by $\Omega_k(z)$, e.g. $\Omega_{(1,1)}(z)=
\frac{dt}{t-z}\wedge\frac{ds}{s-t}$. The projection to $\g$ of this
form will be denoted by $\Omega_k^{\g}(z)$.

The proofs in Section \ref{g} can be modified to get PBW expansions of
$\Omega^{\g}_k(z)$. The only change in the PBW-coefficient results is
that the * of the diagrams has to be decorated by $z$ (instead of
0). E.g. for $k=(2)$, instead of $\asym_{(2)}(dt_1/t_1\wedge
d(t_2-t_1)/(t_2-t_1))=1/(t_1t_2)dt_1\wedge dt_2$ we have
$\asym_{(2)}(dt_1/(t_1-z)\wedge
d(t_2-t_1)/(t_2-t_1))=1/((t_1-z)(t_2-z))dt_1\wedge dt_2$.
For a simple Lie algebra $\g$,
let $V_{\Lambda}$  be a highest weight
$\g$-module with highest weight $\Lambda \in \h^*$ and generating vector
$v_{\Lambda}$.
Recall that the map $U(\n)\to V_{\Lambda}$, $x\to x\cdot v_{\Lambda}$ is surjective.

\begin{definition}
Let $k^{(1)}, k^{(2)},\ldots,k^{(n)}\in \N^r$, $k=\sum
k^{(i)}$. We extend the star multiplication from Section \ref{star}
as follows:
$$*:(\A^{G_{k^{(1)}}}\otimes V_{\Lambda_1}) \otimes
(\A^{G_{k^{(2)}}}\otimes V_{\Lambda_2}) \otimes \ldots \otimes
(\A^{G_{k^{(n)}}}\otimes V_{\Lambda_n}) \to
\A^{G_k}\otimes (V_{\Lambda_1}\otimes V_{\Lambda_2} \otimes \ldots
\otimes V_{\Lambda_n})$$
by
$$(\Omega_1 \otimes v_{1}) *
(\Omega_2 \otimes v_{2}) *
\ldots * (\Omega_n \otimes
v_{n})=(\omega_1 * \ldots * \omega_n) dV_k \otimes (v_1\otimes
v_2\otimes \ldots \otimes v_n),$$
where $\Omega_i=\omega_i dV_{k^{(i)}}$.
\end{definition}

Let $V = V_{\Lambda_1} \otimes \dots \otimes V_{\Lambda_n}$.
We define the $V$-valued differential form of
degree $k$ (c.f.~\cite[(4)]{mukhin}) by

$$
\Omega^V_k=\bigoplus_{k^{(1)}+\ldots+k^{(n)}=k}
\Omega^{\g}_{k^{(1)}}(z_1)v_{\Lambda_1} * \ldots *
\Omega^{\g}_{k^{(n)}}(z_n)v_{\Lambda_n}.
$$

\begin{example} \rm For $n=2$, $r=1$ (i.e. $\g=sl_2$), we have
$$\Omega^V_{(2)}=
\Omega^{\g}_{(2)}(z_1)v_{\Lambda_1} *
\Omega^{\g}_{(0)}(z_2)v_{\Lambda_2} +
\Omega^{\g}_{(1)}(z_1)v_{\Lambda_1} *
\Omega^{\g}_{(1)}(z_2)v_{\Lambda_2} +
\Omega^{\g}_{(0)}(z_1)v_{\Lambda_1} *
\Omega^{\g}_{(2)}(z_2)v_{\Lambda_2} =$$

\begin{align}
& \asym_{(2)}\Big(\frac{dt_1}{t_1-z_1}\wedge \frac{d(t_2-t_1)}{t_2-t_1}\Big)
  \otimes f^2v_{\Lambda_1} \otimes v_{\Lambda_1}
+ \notag \\
& \asym_{(2)}
\Big(\frac{dt_1}{t_1-z_1} \wedge \frac{dt_2}{t_2-z_2}\Big)
 \otimes fv_{\Lambda_1} \otimes
 fv_{\Lambda_2} + \notag \\
& \asym_{(2)}\Big(\frac{dt_1}{t_1-z_2}\wedge
\frac{d(t_2-t_1)}{t_2-t_1}\Big) \otimes v_{\Lambda_1} \otimes f^2v_{\Lambda_2}.
\notag
\end{align}
\end{example}

This can
be visualized by a diagram
$$\Omega^V_{(2)}=
\begin{picture}(20,20)
\put(0,0){$*$}\put(-4,13){${}_{z_1}$}
\put(0,-15){$*$}\put(-4,-18){${}_{z_2}$}
\put(3,3){\line(1,1){15}}
\put(3,3){\line(1,-1){15}}
\put(15,15){$\bullet$}\put(15,24){${}_{1}$}
\put(15,-15){$\bullet$}\put(15,-21){${}_{1}$}
\end{picture}
(f^2v_{\Lambda_1}\otimes v_{\Lambda_2})+
\asym\Big(\ \
\begin{picture}(20,20)
\put(0,13){$*$}\put(0,24){${}_{z_1}$}
\put(3,16){\line(1,0){15}}
\put(15,13){$\bullet$}\put(15,22){${}_{1}$}
\put(0,-12){$*$}\put(0,-16){${}_{z_2}$}
\put(3,-9){\line(1,0){15}}
\put(15,-12){$\bullet$} \put(15,-16){${}_{1}$}
\end{picture}\ \ \Big)
(fv_{\Lambda_1} \otimes fv_{\Lambda_2})
+
\begin{picture}(20,20)
\put(0,0){$*$}\put(-4,-3){${}_{z_2}$}
\put(0,15){$*$}\put(-4,25){${}_{z_1}$}
\put(3,3){\line(1,1){15}}
\put(3,3){\line(1,-1){15}}
\put(15,15){$\bullet$}\put(15,24){${}_{1}$}
\put(15,-15){$\bullet$}\put(15,-21){${}_{1}$}
\end{picture}
(v_{\Lambda_1}\otimes f^2v_{\Lambda_2}),
$$
\ \vskip 0 true cm
\noindent where we also used the PBW expansions from Section \ref{pbwa}.

\end{document}